\documentstyle[12pt,leqno,amssymb]{amsart}

\newtheorem{thm}{Theorem}[section]
\newtheorem{lemma}[thm]{Lemma}
\newtheorem{prop}[thm]{Proposition}

\theoremstyle{definition}
\newtheorem{dfn}[thm]{Definition}
\theoremstyle{remark}

\begin{document}

\newcommand{\xy}{^{x,y}}
\newcommand{\tc}{{\mathfrak{T}}}
\newcommand{\pr}{\protect\ref}
\newcommand{\su}{\subseteq}

\newcommand{\p}{{\overline{p}}}
\newcommand{\C}{{\mathcal C}}
\newcommand{\F}{{\mathcal F}}

\newcommand{\A}{{\mathcal A}}

\newcommand{\R}{{\Bbb R}}
\newcommand{\B}{{\Bbb B}}
\newcommand{\E}{{{\Bbb R}^3}}
\newcommand{\G}{{\Bbb G}}
\newcommand{\Z}{{\Bbb Z}}

\newcommand{\im}{{Imm(F,\E)}}

\newcounter{numb}

\title[Higher Order Invariants of Immersions]
{Higher Order Invariants of Immersions \\ of Surfaces into 3-Space}
\author{Tahl Nowik}
\address{Department of Mathematics, Bar-Ilan University, 
Ramat-Gan 52900, Israel.}
\email{tahl@@math.biu.ac.il}
\date{July 10, 2002}
\thanks{Partially supported by the Minerva Foundation}

\begin{abstract}
We classify all finite order invariants of immersions of a closed
orientable surface into $\E$, with values in any Abelian group. 
We show that they are all functions 
of order one invariants.
\end{abstract}

\maketitle

\section{Introduction}\label{int}

Finite order invariants of stable immersions of a closed orientable
surface into $\E$ have been defined in [N], where all order 1
invariants have been classified. In the present work we classify 
all finite order invariants of order $n>1$, and show that they are 
all functions of the universal order 1 invariant constructed in [N].

The structure of the paper is as follows:
In Section \pr{bac} we summarize the necessary background. 
We define finite order invariants of immersions of a closed orientable
surface into $\E$.
For given surface $F$, regular homotopy class
$\A$ of immersions of $F$ into $\E$, and Abelian group $\G$, we 
define $V_n$ to be the group of all invariants on $\A$
of order at most $n$ with values in $\G$. 
We present  a group $\Delta_n = \Delta_n(\G)$
and an injection $u:V_n/V_{n-1} \to \Delta_n$. The question
of classifying all finite order invariants then becomes the question
of finding the image of $u$. In Section \pr{sta} we state 
our classification. We specify a subgroup $E_n \su \Delta_n$ 
which we claim to be the image of $u$.
In Section \pr{sup} we show that $u(V_n)\supseteq E_n$ by explicitly
constructing a universal finite order invariant.
In Section \pr{sub} we show that $u(V_n)\su E_n$.
The proof relies on the result of Section \pr{sup}.

\section{Background}\label{bac}

In this section we summarize the background needed for this work, 
all of which may be found in [N].
Given a closed orientable surface $F$, $\im$ denotes the 
space of all immersions of $F$ into $\E$, with the $C^1$ topology.
A CE point of an immersion $i:F\to \E$ is a point of self intersection
of $i$ for which the local stratum in $\im$ corresponding to the 
self intersection, has codimension one.
We distinguish twelve types of CEs which we name
$E^0, E^1, E^2, H^1, H^2, T^0, T^1, T^2, T^3, Q^2, Q^3, Q^4$. 
This set of twelve symbols is denoted $\C$.
A co-orientation for a CE is a choice of one of the two sides 
of the local stratum corresponding to the CE.
All but two of the above CE types are nonsymmetric in the
sense that the two sides of the local stratum may be distinguished via
the local configuration of the CE, and for those ten CE types,
permanent co-orientations for the corresponding strata are chosen once and for all. 
The two exceptions are $H^1$ and $Q^2$ which are completely symmetric.

We fix a closed orientable surface $F$ and a regular homotopy class $\A$ of 
immersions of $F$ into $\E$, and
denote by $I_n\su \A$ ($n\geq 0$) the space of all immersions in $\A$ which have precisely 
$n$ CE points (the self intersection being elsewhere stable).
In particular, $I_0$ is the space of all stable immersions in $\A$. 

For an immersion $i:F\to\E$ having a CE located at $p\in\E$, the degree 
$d_p(i) \in \Z$ of $i$ at $p$ is defined in [N]. Then
$C_p(i)$ is the expression $R^a_m$ where $R^a\in\C$
is the symbol describing the configuration of the CE of $i$ at $p$ (one of the twelve symbols) 
and $m=d_p(i)$. 
$\C_n$ denotes the set of all
\emph{un-ordered} $n$-tuples of expressions $R^a_m$ with $R^a\in\C, m\in\Z$.
(So $\C_n$ is the set of un-ordered $n$-tuples of elements of $\C_1$.)
A map $C:I_n \to \C_n$ is defined by
$C(i)=[C_{p_1}(i),\dots, C_{p_n}(i)] \in \C_n$ where $p_1,\dots,p_n$ are the $n$
CE points of $i$. The map 
$C:I_n \to \C_n$ is surjective.

Let $\G$ be any Abelian group and let $f:I_0\to\G$ be an invariant, i.e. a function which is 
constant on each connected component of $I_0$.
Given an immersion $i\in I_n$, a \emph{temporary co-orientation} for $i$ is a choice 
of co-orientation at each of the $n$ CE points $p_1, \dots , p_n$ of $i$.
Given a temporary co-orientation $\tc$ for $i$, and a subset $A\su \{p_1,\dots,p_n\}$,
$i_{\tc,A} \in I_0$ is the immersion obtained from $i$ by resolving all CEs
of $i$ as follows: The CEs at points of $A$ are resolved into the 
positive side with respect to $\tc$,
and those not in $A$ into the negative side.
Given $i\in I_n$ and a temporary co-orientation $\tc$ for $i$,
$f^\tc(i)$ is defined as follows:
$$f^\tc(i)=\sum_{ A \su \{p_1,\dots,p_n\} } (-1)^{n-|A|} f(i_{\tc,A})$$
where $|A|$ is the number of elements in $A$.
The statement $f^\tc(i)=0$ is independent of the temporary co-orientation $\tc$
so we
simply write $f(i)=0$.
An invariant $f:I_0\to\G$ is called \emph{of finite order} if 
there is an $n$ such that $f(i)=0$ for all $i\in I_{n+1}$.
The minimal such $n$ is called the \emph{order} of $f$.
The group of all invariants on $I_0$ of order at most $n$ is denoted $V_n$.

Let $f \in V_n$. If $i\in I_n$ has at least one CE
of type $H^1$ or $Q^2$ and $\tc$ is a temporary co-orientation for $i$,
then $2f^\tc(i)=0$,
and so $f^\tc(i)$ is independent of $\tc$. 
This fact is used to extend any $f\in V_n$ to $I_n$ as $f^\tc(i)$,
where if $i$ includes at least one CE of type $H^1$ or $Q^2$ then
$\tc$ is arbitrary, and if all CEs of $i$ are not of type $H^1$ or $Q^2$
then the permanent co-orientation is used for all CEs of $i$.
(If $f\in V_n$ then we are not extending $f$ to $I_k$ for
$0<k<n$). 
For $f\in V_n$ and $i,j\in I_n$, 
if $C(i)=C(j)$ then $f(i)=f(j)$,
so any $f\in V_n$ induces a well defined
function $u(f):\C_n\to\G$. 
The map $f\mapsto u(f)$ induces an injection $u:V_n / V_{n-1} \to \C_n^*$ 
where $\C_n^*$ is the group of all functions from $\C_n$ to $\G$.
Finding the image of $u$ for all $n$ gives a classification of all finite
order invariants, which is what we do in this work (Theorem \pr{main}).
For order 1 invariants this has been done in [N].

A subgroup $\Delta_n = \Delta_n(\G) \su \C_n^*$ which contains 
the image of $u$ 
is defined as the set of functions in $\C_n^*$ satisfying 
relations which
we write as relations on the symbols $R^a_m$, e.g.
$T^0_m = T^3_m$ will stand for the set of all relations of the form
$g([T^0_m, {R_2}^{a_2}_{d_2},\dots,{R_n}^{a_n}_{d_n}]) =
g([T^3_m, {R_2}^{a_2}_{d_2},\dots,{R_n}^{a_n}_{d_n}])$ with
arbitrary ${R_2}^{a_2}_{d_2},\dots,{R_n}^{a_n}_{d_n} \in \C_1$.
 The relations defining $\Delta_n$
are:
\begin{itemize}
\item $E^2_m = - E^0_m = H^2_m$, \ \  $E^1_m = H^1_m$.
\item $T^0_m = T^3_m$, \ \  $T^1_m = T^2_m$.
\item $2H^1_m =0 $, \ \ $H^1_m = H^1_{m-1}$.
\item $2Q^2_m =0 $, \ \ $Q^2_m = Q^2_{m-1}$.
\item $H^2_m - H^2_{m-1} = T^3_m - T^2_m$.
\item $Q^4_m - Q^3_m = T^3_m - T^3_{m-1}$, \ \  $Q^3_m - Q^2_m = T^2_m - T^2_{m-1}$. 
\end{itemize}

Let $\B\su\G$ be the subgroup defined by
$\B=\{ x\in \G : 2x=0\}$.
To obtain a function $g\in\Delta_1$ 
one may assign arbitrary values in $\G$ for the symbols
$\{T^a_m\}_{a=2,3,m\in\Z}$, $H^2_0$ 
and arbitrary values in $\B$ for the two symbols $H^1_0 , Q^2_0$. Once this is done then 
the value of $g$ on all other symbols is uniquely determined, namely: 
\begin{enumerate}
\item $H^1_m = H^1_0$ for all $m$.
\item $H^2_m = H^2_0 + \sum_{k=1}^m (T^3_k - T^2_k)$ for $m\geq 0$.
\item $H^2_m = H^2_0 - \sum_{k=m+1}^0 (T^3_k - T^2_k)$  for $m<0$.
\item $E^0_m = -H^2_m$, \ \ $E^1_m = H^1_m$, \ \  $E^2_m = H^2_m$ for all $m$.
\item $T^0_m = T^3_m$, \ \ $T^1_m = T^2_m$ for all $m$.
\item $Q^2_m = Q^2_0$ for all $m$.
\item $Q^3_m = Q^2_m + T^2_m - T^2_{m-1}$ for all $m$.
\item $Q^4_m = Q^3_m + T^3_m - T^3_{m-1}$ for all $m$.
\end{enumerate}

A ``universal'' Abelian group $\G_U$ is defined by 
$$\G_U = \left< \{t^a_m\}_{a=2,3, m\in{\Bbb Z}}, h^2_0, h^1_0, q^2_0 \ | \ 
2h^1_0 = 2q^2_0 = 0 \right>.$$ 
Then the universal element $g^U_1\in\Delta_1(\G_U)$ is defined by 
$g^U_1(T^a_m) = t^a_m$ ($a=2,3$),
$g^U_1(H^2_0)=h^2_0$, $g^U_1(H^1_0)=h^1_0$, $g^U_1(Q^2_0)=q^2_0$ 
and the value of $g^U_1$ on all other symbols of $\C_1$ is determined by
formulae 1-8 above. The main result of [N] is that
for any closed orientable surface $F$, regular homotopy class $\A$ of immersions of $F$ 
into $\E$ and Abelian group $\G$, the injection
$u:V_1 / V_0 \to \Delta_1$ is surjective. 
This is shown by constructing
an order one invariant $f^U_1:I_0\to\G_U$ with $u(f^U_1)=g^U_1$. 
Then for arbitrary $\G$, if $g\in \Delta_1(\G)$ then $g=\phi \circ g^U_1$
for $\phi\in Hom(\G_U , \G)$, and $f=\phi \circ f^U_1$ satisfies $u(f)=g$.

\section{Statement of Classification}\label{sta}

Let $X=\{T^a_m\}_{a=2,3,m\in\Z} \cup \{H^2_0 \}$,
and $Y= X \cup \{H^1_0,Q^2_0 \}$
so $X \su Y \su \C_1$.
We observe that as for $\Delta_1$, it is true for any $\Delta_n$
that a function $g\in \Delta_n$ 
may be assigned arbitrary values in $\G$ for 
any un-ordered $n$-tuple of elements of $X$
and arbitrary values in $\B$ for all 
$n$-tuples of elements of $Y$
which include $H^1_0$ or $Q^2_0$ at least once.
Once this is done, 
the value of $g$ on all other $n$-tuples in $\C_n$ is uniquely determined
by the procedure 1-8 of Section \pr{bac}, applied $n$ times, 
independently to each of the $n$ entries of $g$.
This is an evident fact that is somewhat obscured by the fact that we are 
dealing with un-ordered $n$-tuples.
So, in this paragraph only, we revert to ordered notation.
Let us index 
the elements of $Y$
with index $i$ and the elements of $\C_1 \supseteq Y$ with index $j$.
Let us denote $g([i_1,\dots,i_n])$ by $G_{i_1,\dots,i_n}$ where now
the $n$ indices are ordered, and so the expression $G_{i_1,\dots,i_n}$
is symmetric i.e. unchanged when permuting the indices.
The function of $G_{i_1,\dots,i_n}$
produced by the procedure given by 1-8 of 
Section \pr{bac} applied to one of the indices of $G_{i_1,\dots,i_n}$,
is linear with integer coefficients say $y^i_j$ 
where for each $j$, $y^i_j \neq 0$ for only finitely many values of $i$.
Application of the procedure on all $n$ indices of $G_{i_1,\dots,i_n}$
is given by the sum:
$$F_{j_1,\dots,j_n} = \sum_{i_1,\dots,i_n} 
y^{i_1}_{j_1} y^{i_2}_{j_2} \cdots y^{i_n}_{j_n} G_{i_1,\dots,i_n}.$$
With this notation it is clear that $F_{j_1,\dots,j_n}$ thus obtained
is a well defined symmetric function of $j_1,\dots,j_n$ which is the 
unique extension of $G_{i_1,\dots,i_n}$ 
to the indices $j_1,\dots,j_n$, satisfying the relations defining $\Delta_n$.

We will now define $E_n \su \Delta_n$ by two additional
restrictions on the functions $g\in \Delta_n$.
Thanks to the discussion of the previous paragraph, we
may state the additional restrictions in terms of the
values of $g$ on $n$-tuples of elements of $Y$ only.
Given an un-ordered $n$-tuple $z$ of elements of $Y$, 
we define $m_{H^1_0}(z)$ and $m_{Q^2_0}(z)$ 
as the number of times that
$H^1_0$ and $Q^2_0$ appear in $z$ respectively.
We define $r(z)$, (the \emph{repetition} of $H^1_0$ and $Q^2_0$ in $z$), as 
$$r(z)= \max(0, m_{H^1_0}(z)-1) + \max(0, m_{Q^2_0}(z)-1).$$

\begin{dfn}\label{En}
Given an Abelian group $\G$,
$E_n = E_n(\G) \su \Delta_n(\G)$ is the subgroup consisting
of all $g\in \Delta_n(\G)$
satisfying the following two additional restrictions:

\begin{enumerate}

\item When $n \geq 3$, $g$ must satisfy the relation $H^1_0 H^1_0 Q^2_0 = H^1_0 Q^2_0 Q^2_0$.
\newline By this we mean that 
$g([H^1_0, H^1_0, Q^2_0, {R_4}^{a_4}_{d_4},\dots,{R_n}^{a_n}_{d_n}]) =
g([H^1_0, Q^2_0, Q^2_0, {R_4}^{a_4}_{d_4},\dots,{R_n}^{a_n}_{d_n}])$ for
arbitrary ${R_4}^{a_4}_{d_4},\dots,{R_n}^{a_n}_{d_n} \in Y$.
\label{e1}

\item For any un-ordered $n$-tuple $z$ of elements of $Y$,
$g(z)\in 2^{r(z)}\G$, i.e. there exists an element $a\in \G$ such
that $g(z)=2^{r(z)} a$.
(Note that whenever $r(z) > 0$ then in particular $H^1_0$ or $Q^2_0$ does appear
in $z$ so in fact we have $g(z) \in \B \cap 2^{r(z)}\G$.)
\label{e2}

\end{enumerate}
\end{dfn}

In this work we prove:

\begin{thm}\label{main}
For any closed orientable surface $F$, regular homotopy class $\A$ of immersions of $F$ 
into $\E$ and Abelian group $\G$, the 
image of the injection $u:V_n / V_{n-1} \to \Delta_n$ is $E_n$.

Furthermore, for any $g\in E_n$ there exists a function (not homomorphism)
$s: \G_U \to \G$ such that the invariant $f=s\circ f^U_1$ is of order $n$ and
satisfies $u(f)=g$. ($f^U_1$ of Section \pr{bac}.)
It follows that all finite order invariants are functions
of the order 1 invariant $f^U_1$.
\end{thm}

\section{Proof that $u(V_n) \supseteq E_n$}\label{sup}

For convenience, we rename the generating elements of the group
$\G_U$. 
\newline 
We relabel 
$\{t^a_m\}_{a=2,3, m\in{\Bbb Z}} \cup \{ h^2_0 \}$ as 
$\{ a_i \}_{i\in X}$ where $X$ is a countable set of indices,
and relabel $h^1_0, q^2_0$ as $b,c$.
We define algebraic structures $K \su L \su M$, where $L$ is a commutative 
ring, $K$ is a subring of $L$, and $M$ is a module over $K$.
$L$ is defined as the ring of formal power series with integer coefficients
and variables
$\{ a_i \}_{i\in X} \cup \{b,c\}$ and with relations 
\begin{itemize}
\item $b^2 c = b c^2$. 
\item $2b=2c=0$.
\end{itemize}
We emphasize that though there is an infinite
set of variables, any given power series may include only finitely many
monomials of any given degree $n$.
Given a monomial $p$, 
we define $m_b(p)$ and $m_c(p)$ as the multiplicity of $b$ and $c$ in $p$ respectively. We 
define $r(p)$, (the repetition of $b$ and $c$ in $p$), as 
$$r(p)= \max(0, m_b(p)-1) + \max(0, m_c(p)-1).$$ 
$r(p)$ is preserved under the 
relations in $L$ and so is well defined on 
equivalence classes of monomials.
The equivalence class of a monomial $p$ 
will be denoted $\p$. 
We note that an equivalence class
$\p$ includes more than just the one monomial $p$ iff
$m_b(p)\geq 1$, $m_c(p)\geq 1$ and $r(p)\geq 1$;
it then includes precisely $r(p)+1$ different monomials.
When $\p$ includes only the one monomial $p$ then we will
interchangeably write $\p$ and $p$.
Now $K\su L$ is defined to be the subring of power series
including only the variables $\{ a_i \}_{i\in X}$.
On the other hand 
we extend $L$ to a larger structure $M$ which will 
be a module over the subring $K$, 
as follows:
For each $\p$ for which $p$ is a monomial with coefficient 1,
we adjoin a new element $\zeta_\p$ satisfying the
relation $2^{r(\p)} \zeta_\p = \p$.  
The new elements $\zeta_\p$ will be considered monomials of the same degree as
$\p$, and will appear as terms in our formal power series.
(Indeed, one can think of $\zeta_\p$ as $2^{-r(\p)}\p$.)
Note that if $r(p)=0$ then $\zeta_p=p$,
in particular, $\zeta_1 = 1$ and $\zeta_e=e$ for each generating variable $e$.
Now $K$ acts on $M$ as follows: If $k\in K, p\in L$ are monomials
then $k \cdot \zeta_\p = \zeta_{\overline{kp}}$. This is extended in the natural way
to an action of power series in $K$ on power series in $M$. 
We note that the whole of $L$ cannot act on $M$ in this way, since we would get
contradictions such as $0\neq b^2 = 2 \zeta_{b^2} = 2b \cdot \zeta_b =0$.
In particular, we do not have a ring structure on $M$.
For each $n\geq 0$ we denote by $K_n \su L_n \su M_n$ the additive subgroups
of $K\su L \su M$ respectively generated by the monomials of degree $n$.
(As mentioned, the degree of $\zeta_\p \in M$ is defined as the degree of $\p$).
We note $L_1=M_1$ is the Abelian group with generators 
$\{ a_i \}_{i\in X} \cup \{b,c\}$ and relations $2b=2c=0$,
i.e. the group $\G_U$.
We have $L_1=K_1\oplus S$ where $S \su L_1$
is the four element subgroup generated by $b,c$.
We now define a function
$\F:L_1\to M$ as follows: We first define $\F:K_1 \to K$ as the group 
homomorphism
from the additive group $K_1$ to the multiplicative group of 
invertible elements in $K$, which is given on generators
by $$\F(a_i)=\sum_{n=0}^\infty a_i^n.$$ These are indeed invertible elements, giving
$$\F(-a_i)=\bigg(\sum_{n=0}^\infty a_i^n\bigg)^{-1} = 1-a_i.$$
Note that for $x=\pm a_i$, $\F(x)=1+x+T_2$ where $T_2$ stands for 
the ``higher order terms'' (or ``tail'') of the given series,
i.e. a power series including only monomials of degree at least 2.
It follows that for any $x\in K_1$,
$\F(x)=1+x+T_2$.

We then define $\F:S\to M$ explicitly on the four elements of $S$ as follows: 
\begin{enumerate}
\item $\F(0)=1$.
\item $\F(b)=\sum_{n=0}^\infty \zeta_{b^n}$.
\item $\F(c)=\sum_{n=0}^\infty \zeta_{c^n}$.
\item $\F(b+c)=1+b+c+\sum_{n=2}^\infty (\zeta_{b^n} + \zeta_{c^n} + \zeta_{\overline{bc^{n-1}}})$.
\end{enumerate}

Finally, $\F:L_1 \to M$ is defined as follows: Any element in $L_1$ is uniquely written as
$k+s$ with $k\in K_1, s\in S$, and we define $\F(k+s)=\F(k)\F(s)$
where the product on the right is the action of $K$ on $M$. 
It follows that for any $k\in K_1, l\in L_1$: 
$\F(k+l)=\F(k)\F(l)$.

For any $(n+1)$-tuple $(l; \ l_1,l_2,\dots,l_n)$ of elements of $L_1$, we define
$$\F'(l; \ l_1,\dots,l_n)=
\sum_{A\su \{1,\dots,n\}} (-1)^{n-|A|} \F(l+\sum_{i\in A} l_i).$$

\begin{prop}\label{al}
For any $n+1$-tuple $(l; \ l_1,\dots,l_n)$ of elements of $L_1$:  
$$\F'(l; \ l_1,\dots,l_n)= l_1 l_2 \cdots l_n + T_{n+1}$$
where $T_{n+1}$ stands for higher order terms
i.e. some power series in $M$ having only monomials of degree at least $n+1$.
\end{prop}

\begin{pf}
We first show that it is enough to prove the proposition
for the case when all $l_i$ are generators or minus generators
of $L_1$, i.e. of the form $\pm e$ where $e\in \{a_i\}_{i\in X} \cup \{b,c \}$.
We prove this by induction on the sum of the lengths of $l_1,\dots,l_n$
in terms of the generators. 
Say $l_1$ is not a generator and 
$l_1=l'_1 + l_1''$ where $l'_1$ and $l_1''$ have shorter
length than $l_1$ in terms of the generators. Then
\begin{align*}
\F'(l; \ l_1,\dots,l_n)= &
\sum_{A\su \{1,\dots,n\}} (-1)^{n-|A|} \F(l+\sum_{i\in A} l_i) \\ = &
\sum_{A\su \{2,\dots,n\}} (-1)^{n-|A|-1} 
\bigg(\F(l+l_1+\sum_{i\in A} l_i) - \F(l+\sum_{i\in A} l_i)\bigg) \\ = & 
\sum_{A\su \{2,\dots,n\}} (-1)^{n-|A|-1} 
\bigg(\F(l+l'_1+l_1''+\sum_{i\in A} l_i) - \F(l+l'_1+\sum_{i\in A} l_i)\bigg) + \\
 & \qquad \qquad + \sum_{A\su \{2,\dots,n\}} (-1)^{n-|A|-1} 
\bigg(\F(l+l'_1+\sum_{i\in A} l_i) - \F(l+\sum_{i\in A} l_i)\bigg) \\ = & 
\F'(l+l'_1; \ l_1'',l_2,\dots,l_n) + \F'(l; \ l'_1,l_2,\dots,l_n) \\ = &
l_1'' l_2 \cdots l_n +T_{n+1} \ \ + \ \  l'_1 l_2 \cdots l_n +T_{n+1}  \ =
\ l_1 l_2 \cdots l_n +T_{n+1}.
\end{align*}

We are left with proving the proposition in case all $l_i$ are generators
or minus generators. We prove this by induction on $n$.
If one of the $l_i$,
say $l_1$,
is in $K_1$, then 
\begin{align*}
\F'(l; \ l_1,\dots,l_n)= &
\sum_{A\su \{1,\dots,n\}} (-1)^{n-|A|} \F(l+\sum_{i\in A} l_i) \\ = &
\sum_{A\su \{2,\dots,n\}} (-1)^{n-|A|-1} 
\bigg(\F(l+l_1+\sum_{i\in A} l_i) - \F(l+\sum_{i\in A} l_i)\bigg) \\ = &
(\F(l_1)-1)\sum_{A\su \{2,\dots,n\}} (-1)^{n-1-|A|} 
\F(l+\sum_{i\in A} l_i) \\ = &
(\F(l_1)-1)\F'(l; \ l_2,\dots,l_n)  \\ = &
(l_1 + T_2)(l_2\cdots l_n + T_n)  = 
l_1 l_2 \cdots l_n + T_{n+1}.
\end{align*}
Note that the term multiplying on the left is indeed an element of $K$.

So assume now that all $l_1,\dots,l_n$ are $b$ and $c$. Assume $k$ of them are $b$
and $n-k$ of them are $c$. We first deal with the case $k=n$ i.e. $l_1,\dots,l_n$ are all $b$ 
(the case $k=0$ is identical). Since $2b=0$ we get:
\begin{align*}
\F'(l; \ b,b,\dots,b)= &
\sum_{A\su \{1,\dots,n\}} (-1)^{n-|A|} \F(l+|A|b) \\ = &
\pm \bigg(\sum_{|A| \text{odd}} \F(l+b) - \sum_{|A| \text{even}} \F(l)\bigg)  = 
\pm 2^{n-1} \bigg( \F(l+b) - \F(l) \bigg).
\end{align*}
since indeed there is an equal number of odd and even sized subsets 
of $\{1,\dots,n\}$ (e.g. since $\sum_{i=0}^n (-1)^i {n \choose i} = 0$).

Now letting
$l=k+s$ with $k\in K_1$, $s\in S$, 
we get: 
\begin{align*}
 \pm 2^{n-1} (\F(k+s+b) - \F(k+s))  = & \pm\F(k)2^{n-1} (\F(s+b) - \F(s)) \\
= & \pm (1+T_1) 2^{n-1} (\F(s+b)-\F(s)).
\end{align*}
Since multiplication by $1+T_1 \in K$ leaves the lowest order term unchanged
we may assume we have only $\pm 2^{n-1} (\F(s+b) -\F(s))$. 
If $s=0$ or $b$ then 
$$\pm 2^{n-1} (\F(s+b) -\F(s)) = \pm 2^{n-1} \sum_{m=1}^\infty \zeta_{b^m}
= b^n + T_{n+1}$$ since $r(b^m)=m-1$ and so $2^{n-1}\zeta_{b^m}=0$ for $m<n$ and
$2^{n-1}\zeta_{b^n}=b^n$. (The $\pm$ was dropped since $2b^n=0$.)
If $s=c$ or $c+b$ then we get 
$$\pm 2^{n-1} (\F(s+b) -\F(s)) = \pm 2^{n-1} \bigg( b + 
\sum_{m=2}^\infty(\zeta_{b^m} + \zeta_{\overline{bc^{m-1}}})\bigg)
= b^n + T_{n+1}$$
since again $r(b^m)=m-1$,
but also $r(bc^{m-1})=m-2$ and so $2^{n-1} \zeta_{\overline{bc^{m-1}}} =0$ for $m\leq n$.

We are left with the case of $b$ appearing $k$ time and $c$ appearing $n-k$ times, 
with $0<k<n$. Since $2b=2c=0$, we get:
\begin{align*}
\F'(l; \ b,\dots,b,c,\dots,c)= & 
\sum_{B\su \{1,\dots,k\},C\su \{k+1,\dots,n\}} (-1)^{n-|B|-|C|} \F(l+|B|b+|C|c) \\ = &
\pm\bigg(
\sum_{|B| \text{odd}, |C| \text{odd}} \F(l+b+c) 
-\sum_{|B| \text{odd}, |C| \text{even}} \F(l+b) \\
& \qquad \qquad -\sum_{|B| \text{even}, |C| \text{odd}} \F(l+c) 
+\sum_{|B| \text{even}, |C| \text{even}} \F(l) 
\bigg) \\ = &
\pm 2^{n-2}\bigg(\F(l+b+c)-\F(l+b)-\F(l+c)+\F(l)\bigg).
\end{align*}
As before we may assume $l=s\in S$. For each of the four elements
$s\in S$ we get:
\begin{align*}
\pm 2^{n-2}\bigg(\F(s+b+c)-\F(s+b)-\F(s+c)+\F(s)\bigg) = &
\pm 2^{n-2} \sum_{m=2}^\infty \zeta_{\overline{bc^{m-1}}}  \\ = &
\overline{bc^{n-1}} + T_{n+1} = \overline{b^k c^{n-k}} + T_{n+1} 
\end{align*}
since $r(\overline{bc^{m-1}})=m-2$ and so 
$2^{n-2}\zeta_{\overline{bc^{m-1}}}=0$ for $m<n$,
and $2^{n-2}\zeta_{\overline{bc^{n-1}}}= \overline{bc^{n-1}}$. 

\end{pf}

We now return to our original symbols $\{t^a_m\}_{a=2,3, m\in{\Bbb Z}}, h^2_0$ and 
$h^1_0, q^2_0$ in place of $\{a_i \}_{i\in X},b,c$.
Extending the definition of $g^U_1$ of Section \pr{bac}
($g^U_1: \C_1 \to \G_U = L_1$),
let $g^U_n \in E_n(M_n)$ be first defined 
on un-ordered $n$-tuples of elements of 
$Y$
by
$g^U_n([{R_1}^{a_1}_{d_1},\dots,{R_n}^{a_n}_{d_n}])=
\overline{{r_1}^{a_1}_{d_1} \cdots {r_n}^{a_n}_{d_n}} \in M_n$
where $r_i$ is the lower case letter corresponding to the capital letter $R_i$.
(Note that the upper indices are not powers but part of the given symbols.)
The value of $g^U_n$ on all other $n$-tuples in $\C_n$ is then determined by
the procedure discussed in the opening paragraph of Section \pr{sta}.
It follows that for \emph{any} $n$-tuple $[{R_1}^{a_1}_{d_1},\dots,{R_n}^{a_n}_{d_n}]\in\C_n$,
$g^U_n([{R_1}^{a_1}_{d_1},\dots,{R_n}^{a_n}_{d_n}])
=g^U_1({R_1}^{a_1}_{d_1})g^U_1({R_2}^{a_2}_{d_2})\cdots g^U_1({R_n}^{a_n}_{d_n})$
where the product is that in $L$. By construction of $M_n$ indeed $g^U_n \in E_n(M_n)$. 

Let $\F_n:L_1 \to M_n$ be the projection onto $M_n$ of $\F:L_1 \to M$,
and let $f^U_n: I_0 \to M_n$ 
be the invariant given by $f^U_n = \F_n \circ f^U_1$.
($f^U_1 : I_0 \to \G_U = L_1$ appearing in 
the concluding paragraph of Section \pr{bac}.)
Let $i\in I_m$, where $m$ is either $n$ or $n+1$, with
CEs at $\{ p_1,\dots,p_m \}$
and let $\tc$ be a temporary co-orientation for $i$ which
is the \emph{permanent} co-orientation 
wherever there is one. Then by  
Proposition \pr{al} and since $u(f^U_1)=g^U_1$:
\begin{align*}
f^U_n(i)=(f^U_n)^\tc(i)=&\sum_{ A \su \{p_1,\dots,p_m\} } (-1)^{m-|A|} \F_n\circ f^U_1(i_{\tc,A}) \\ = &
\sum_{ A \su \{p_1,\dots,p_m\} } (-1)^{m-|A|} 
\left( \F_n \bigg( f^U_1(i_\varnothing) + \sum_{p_j \in A} g^U_1(C_{p_j}(i))
  \bigg) \right) \\ = &
\begin{cases} 
0 & \text{ \ when \ } m=n+1 \\
g^U_1(C_{p_1}(i)) g^U_1(C_{p_2}(i))  \cdots g^U_1(C_{p_n}(i)) 
= g^U_n(C(i)) & \text{ \ when \ } m=n.
\end{cases}
\end{align*}
That is $f^U_n$
is an invariant of order $n$ with $u(f^U_n)=g^U_n$.
Now for arbitrary Abelian group $\G$, if $g\in E_n(\G)$ then there is
$\phi \in Hom(M_n,\G)$ such that $g=\phi \circ g^U_n$.
It is then clear that $f=\phi \circ f^U_n$ ($=\phi \circ \F_n \circ f^U_1$) is an invariant
of order $n$ with $u(f)=g$.
This proves that $u(V_n)\supseteq E_n$ for any $\G$.
To complete the proof of Theorem \pr{main} it remains to show that
$u(V_n)\su E_n$. When this is established 
then $\F \circ f^U_1$ may be named a ``universal finite order invariant''.

\section{Proof that $u(V_n) \su E_n$}\label{sub}

For $x,y,n \geq 0$
we define $I\xy_n$ to be the space of immersions in $I_{x+y+n}$ with $x$
designated CEs of type $H^1_0$ with choice of ordering on them,
$y$ designated CEs of type $Q^2_0$ with choice of ordering on them, 
and a choice of co-orientation for these $x+y$ CEs. The remaining $n$ CEs may be of any type and they are not
ordered nor co-oriented. 
So the same underlying immersion appears $2^{x+y} x! y!$ times
in $I\xy_n$ with different
choices of ordering and co-orientations. Also note that $I^{0,0}_n=I_n$.

An $(x,y)$-invariant is a function $f:I\xy_0 \to \G$ which is constant
on the connected components of $I\xy_0$. 
Given a temporary co-orientation $\tc$ for the $n$ non-designated CEs of 
$i\in I\xy_n$,
we define $i_{\tc,A} \in I\xy_0$  as before,
resolving only non-designated CEs and 
keeping the order and co-orientation of the designated CEs.
We may then define $f^\tc(i)$ and 
invariants of order $n$ as before,
and define $V\xy_n$ to be the space of all $(x,y)$-invariants of order at most $n$.
We define $C:I\xy_n \to \C_n$ as before, using the $n$ non-designated CEs.
Again $C$ is surjective and
induces an injection $u: V\xy_n / V\xy_{n-1} \to \Delta_n$. Indeed all arguments
(appearing in [N]) showing that $u$ may be defined on $V_n$ and that
$u(V_n)\su \Delta_n$, are applicable 
in just the same way to show that the same is true for $V\xy_n$.
(As first step note that
by [N], namely Proposition 3.4, proof of Proposition 3.5, and Remark 3.7,
for any $i,j \in I\xy_n$, $C(i)=C(j)$ iff there is an AB equivalence 
between the underlying
immersions which preserves all additional structure, i.e. it brings each designated 
CE of $i$ to its counterpart in $j$, and with the right co-orientation.
The notion of AB equivalence appears in [N] Definition 3.3, which also refers
to Definition 2.1.)
We will show that in fact $u(V\xy_n)\su E_n$ for any $x,y$. In particular we
will have $u(V_n) = u(V^{0,0}_n) \su E_n$ which is the aim of this section.

\begin{dfn}\label{tag}
If $i\in I\xy_0$ then we denote by $i'$ the immersion in $I_0$
which is obtained from $i$ by resolving all $x+y$ designated
CEs into the positive side determined by their chosen co-orientation.
\end{dfn}
 
\begin{lemma}\label{ext}
Given $n$, assume it is known that for any $k<n$, (and any $x,y$)
$u(V\xy_k)\su E_k$. Then for any $k<n$ and any $f\in V\xy_k$ there
exists $F\in V_k$ such that for any $i\in I\xy_0$, $f(i)=F(i')$.
\end{lemma}

\begin{pf}
Proof by induction on $k$ ($< n$). 
By our assumption $u(f) \in E_k$. By Section \pr{sup}, $u(V_k) \supseteq E_k$ and so
there exists $G \in V_k$ with $u(G)=u(f)$. 
Let $h$ be the invariant on
$I\xy_0$ defined by $h(i) = f(i) - G(i')$. 
Then $u(h)=0$ ($u$ defined on $V\xy_k$) so $h \in V\xy_{k-1}$, so by the induction 
hypothesis there 
is $H \in V_{k-1}$ such that $h(i) = H(i')$ for all $i \in I\xy_0$.
$F=H+G$ is the required invariant on $I_0$.
\end{pf}

\begin{lemma}\label{two}
Given $n$, assume it is known that for any $k<n$, (and any $x,y$)
$u(V\xy_k)\su E_k$.
Let $f\in V\xy_n$ and
let $i,j \in I\xy_0$ be two immersions such that there is an AB equivalence between
them (respecting ordering and co-orientations of the designated CEs) 
during which precisely two CEs occur,
both of which are of type $H^1_0$.
Then $f(i)=f(j)$. 

The same is true for $Q^2_0$.
\end{lemma}

\begin{pf}
Given $f\in V\xy_n$ we define $f^H \in V^{x+1,y}_{n-1}$ and $f^Q \in V^{x,y+1}_{n-1}$
as follows: For $i\in I^{x+1,y}_0$ let $f^H(i) = f(i^+) - f(i^-)$ where 
$i^+ \in I\xy_0$ is the immersion obtained from $i$ by resolving the $(x+1)$th
designated CE of type $H^1_0$ into the positive side determined by the chosen 
co-orientation, and the ordering and co-orientation on the remaining designated
CEs remains as in $i$. Similarly $i^-$ is defined using the negative side of
the co-orientation at the same CE. In the same way $f^Q$ is defined 
on $i\in I^{x,y+1}_0$ using the $(y+1)$th
designated CE of type $Q^2_0$. Indeed it is clear that 
$f^H \in V^{x+1,y}_{n-1}$ and $f^Q \in V^{x,y+1}_{n-1}$.
We continue discussing $H^1_0$ but clearly all will be true for $Q^2_0$ as well.
By our assumption and Lemma \pr{ext} there exists $G\in V_{n-1}$
such that $f^H(i) = G(i')$ for all $i\in I^{x+1,y}_0$. 
Let $J_t:F\to\E$ ($0\leq t \leq 1$)
be the AB equivalence in the assumption of the lemma so $J_0=i, J_1=j$ and
assume the two CEs occur at times $1\over 3$ and $2\over 3$.
We make $J_{1\over 3}$ and $J_{2\over 3}$ into elements of $I^{x+1,y}_0$
by announcing the additional CE that is occurring as the $(x+1)$th designated
CE of type $H^1_0$. For $J_{1\over 3}$ we choose the co-orientation
of this $(x+1)$th CE to be represented by the motion of $J_t$ through $J_{1\over 3}$
with increasing time, whereas for 
$J_{2\over 3}$ we use the motion of $J_t$ with \emph{decreasing} time.
So $J_{1\over 2}$ is on the positive side of both $J_{1\over 3}$ and $J_{2\over 3}$.
The co-orientation and order on
all other designated CEs of $J_{1\over 3}$ and $J_{2\over 3}$ are those of
$i$ which are continuously carried along the regular homotopy $J_t$.
Similarly $J_{1\over 2}$ is made into an element of $I\xy_0$ by 
continuously carrying the co-orientation and order of the CEs of $i$
along $J_t$.
We get 
$$f(J_{1\over 2}) - f(i) = 
f^H(J_{1\over 3}) = G(J'_{1\over 3}) = G(J'_{2\over 3}) =
f^H(J_{2\over 3}) =
f(J_{1\over 2}) - f(j)$$
since $G$ is defined on $I_0$
and $J'_{1\over 3}$ and $J'_{2\over 3}$ are in the same connected component
of $I_0$. And so we get $f(i)=f(j)$.
\end{pf}

We now prove that $u(V\xy_n) \su E_n$, by induction on $n$
i.e. we assume it is true for any $k<n$, and so the conclusion
of Lemma \pr{two} holds.
Let $i\in I\xy_n$ having all its non-designated CEs from the set $Y$
and located at $p_1,\dots,p_n$.
If $\tc$ is a temporary co-orientation for $p_1,\dots,p_n$ then 
by Lemma \pr{two}, $f(i_{\tc,A})=f(i_{\tc,B})$ whenever
$A,B \su \{p_1,\dots,p_n\}$
have the same number mod 2 of points of type $H^1_0$ and
the same number mod 2 of points of type $Q^2_0$, all other points being the same.
It follows that in the expression for $f^\tc(i)$, each value of
$(-1)^{n-|A|} f(i_{\tc,A})$ appears $2^r$ times, where $r=r(C(i))$,
and so we get the following: 

\begin{lemma}\label{tr}
For $i$ as above let $R\su \{ p_1,\dots,p_n \}$ be a subset which includes all
the points which are not of type $H^1_0$ and $Q^2_0$, and includes
only one of the points of type $H^1_0$ if such exists, and only one
of the points of type $Q^2_0$ if such exists, then:
$$f^\tc(i)=  \sum_{ A \su \{p_1,\dots,p_n\} } (-1)^{n-|A|} f(i_{\tc,A}) =
 2^r \sum_{A\su R } (-1)^{n-|A|} f(i_{\tc,A}).$$
\end{lemma}

This proves that $u(f)$ satisfies property \pr{e2} of Definition \pr{En}.

Now let $i\in I\xy_{n+1}$ having all its non-designated CEs from the set $Y$
and located at $\{p_1,\dots,p_{n+1}\}$ and assume $C_{p_1}(i)=C_{p_2}(i)=H^1_0$
and $C_{p_n}(i)=C_{p_{n+1}}(i)=Q^2_0$.
Given a temporary co-orientation $\tc$ for $i$,
let $i_0\in I\xy_n$ be the immersion obtained from $i$ by resolving the 
CE at $p_{n+1}$ into the negative side determined by $\tc$.
Let $i_1 \in I\xy_n$ be similarly defined using the point $p_1$.
So the CEs of $i_0$ are $\{p_1,\dots,p_n\}$ and the CEs 
of $i_1$ are $\{p_2,\dots,p_{n+1}\}$. Let $\tc_0, \tc_1$ be the temporary
co-orientations for $i_0, i_1$ respectively, which are the restrictions of 
$\tc$, then for any $A\su \{p_2,\dots,p_n\}$, $(i_0)_{\tc_0,A} = i_{\tc,A} = (i_1)_{\tc_1,A}$.
Let $R\su \{p_2,\dots,p_n \}$ be the set including $p_2,p_n$ and all points
which are not of type $H^1_0$ and $Q^2_0$.
Then this $R$ may be used as the $R$ appearing in Lemma \pr{tr},
for both $i_0$ and $i_1$. Furthermore $r(C(i_0))=r(C(i_1))$ which we denote $r$,
so we get by Lemma \pr{tr}: 
$$
f(i_0)=f^{\tc_0}(i_0) =  
2^r \sum_{A\su R } (-1)^{n-|A|} f(i_{\tc,A}) 
= f^{\tc_1}(i_1)=f(i_1).
$$
This proves that $u(f)$ satisfies property \pr{e1} of Definition \pr{En} since 
$p_1,p_2,p_n$ are CEs of $i_0$
and $p_2,p_n,p_{n+1}$ are CEs of $i_1$,
all other CEs being the same,
and since any equality 
stated by \pr{e1} of Definition \pr{En} may indeed be realized by
an immersion $i\in I\xy_{n+1}$ with such set $p_1,\dots,p_{n+1}$ of CEs.
This completes the proof of Theorem \pr{main}.

\end{document}